\documentclass{article}
\def\ni{\noindent}
\def\Z{\hbox{\bf Z}}
\def\Q{\hbox{\bf Q}}
\def\C{\hbox{\bf C}}
\def\mod{\mathop{\;\rm mod}\nolimits}
\def\Pic{\mathop{\rm Pic}\nolimits}

\def\ms{\mapsto}

\begin{document}

\begin{center}
{\sc \large{Couples de jacobiennes isog\`enes de courbes hyperelliptiques de
genre
arbitraire}}\\
\medskip

{\sc Jean-Fran\c cois Mestre, Universit\'e Paris $7$, Paris.}
\end{center}

\medskip

\section{Introduction}
Soit $C$ une courbe de genre $g$, $J(C)$ sa jacobienne, $H$ un sous-groupe
totalement isotrope de rang $g$ de $J(C)[2]$; la vari\'et\'e ab\'elienne 
quotient $A$ de $J(C)$ par $H$ est principalement polaris\'ee, 
mais, pour $g\geq 4$, n'est en g\'en\'eral pas
une jacobienne. A fortiori, si $C$ est  hyperelliptique, et $g\geq 3$, $A$ n'est en
g\'en\'eral pas la jacobienne d'une courbe hyperelliptique. 

Il ne
semble m\^eme pas connu que, pour $g$ assez grand, il existe   au moins un
couple de courbes hyperelliptiques $(C,C')$ de genre  $g$ dont les
jacobiennes sont reli\'ees par une $2\times \ldots \times 2$ isog\'enie. 

Notons n\'eanmoins que B. Smith a obtenu des familles \`a $3$ (resp. $2$,
resp. $1$) param\`etres
de tels couple de courbes de genre $6,12,14$ (resp. $3,6,7$, resp. $5,10,15$).

\medskip

Nous montrons ici que, pour tout $g\geq 2$, il existe une famille \`a
$g+1$ param\`etres de paires de courbes hyperelliptiques $(C,C')$ dont les
jacobiennes sont reli\'ees par une isog\'enie de noyau isomorphe
\`a $(\Z/2\Z)^g$. Plus pr\'ecis\'ement:

\medskip

\ni
{\sc Th\'eor\`eme .-} {\it Soit $g$ un entier  $\geq 1$, et $K=\Q(a_1,\ldots,a_g,v)$, o\`u  
$a_1,\ldots,a_g,v$ sont des ind\'etermin\'ees; il existe une correspondance $2-2$ entre les courbes $C$ et $C'$ d'\'equations respectives

$$y^2=(x-v)(vx-1)(x^2-a_1)\ldots (x^2-a_g)$$
et 
$$y^2=(x-v)(vx-(-1)^g)(x^2-b_1)\ldots (x^2-b_g),$$
o\`u $b_i=\frac{a_i v^2-1}{a_i-v^2}$ pour $1\leq i\leq g$,
induisant une $\overbrace{2\ldots 2}^g$ isog\'enie  entre leurs jacobiennes. 

\medskip
La jacobienne de $C$ est absolument simple; de plus,
lorsqu'on sp\'ecialise les $a_i$ et $v$ en des \'el\'ements de $\C$, l'image des courbes
$C$ dans la vari\'et\'e des modules des courbes hyperelliptiques de genre $g$ sur $\C$ est de dimension $g+1$.   
}

\bigskip
{\sc Remarques.-} 1) Pour $g$ pair par exemple, ceci permet d'obtenir une famille de dimension $g/2+1$ de courbes hyperelliptiques dont l'anneau des endomorphismes de la jacobienne 
contient $\Z[\sqrt{2}]$: $v$ et $a_i$ ($1\leq i \leq g/2$) \'etant arbitraires, 
il suffit de prendre
$a_{g/2+i}=\frac{a_iv^2-1}{a_i-v^2}$ pour $1\leq i\leq g/2.$

2) Dans le cas du genre $2$, on retrouve la correspondance de
Richelot (cf. par exemple [1], [2], [3]).

\medskip

\section{D\'emonstration}

On garde les notations du th\'eor\`eme. Notons $p_0(x)=q_0(x)=(x-v)(vx-1)$ et, pour $1\leq i\leq g$, $p_i(x)=x^2-a_i$ et
$q_i(x)=x^2-b_i$; si l'on pose $S(x,z)=x^2z^2-v^2(x^2+z^2)+1$, o\`u $z$ est une ind\'etermin\'ee,
on a les 
identit\'es
$$\left\{\begin{array}{l}p_2(v)p_1(x)q_2(z)-p_1(v)p_2(x)q_1(z)+(a_1-a_2)S(x,z)=0\\
(1-v^2)S(x,z)=2p_0(x)q_0(z)-(v^2+1)(1-xv-zv+xz)^2
\end{array}\right.$$

d'o\`u

$$\left\{\begin{array}{ll}
p_2(v)p_1(x)q_2(z)&\equiv p_1(v)p_2(x)q_1(z)\\
2p_0(x)q_0(z)&\equiv (v^2+1)(1-xv-zv+xz)^2\end{array}\right.\mod S.$$

\subsection{Le cas o\`u $g$ est pair}

Supposons d'abord $g$  pair; pour $1\leq i\leq g$,  on a
d'apr\`es la formule pr\'ec\'edente $p_{2i}(v)p_{2i-1}(x)q_{2i}(z)\equiv p_{2i-1}(v)p_{2i}(x)q_{2i-1}(z)\mod  S$
pour $1\leq i\leq g/2$. 

Par suite, si $M(x,z)=p_2(v)p_4(v)\ldots p_{g}(v)p_1(x)q_2(z)p_3(x)q_4(z)\ldots
p_{g-1}(x)q_g(z)$, on a $$\prod_{i=1}^g p_i(v)p_i(x) \prod_{i=1}^g q_i(z)\equiv M(x,z)^2\mod S.$$

Si $C$ est la courbe d'\'equation $y^2=A\prod_{i=0}^g p_i(x)$, o\`u $A=2(v^2+1)\prod_{i=1}^g p_i(v)$,
et $C'$ celle d'\'equation
$t^2=\prod_{i=0}^g q_i(z)$, on a donc une correspondance $\Gamma$ sur $C\times C'$ 
d\'efinie par
les \'equations

$$\Gamma:\;\;\;\left\{\begin{array}{ll}(S(x,z)&=0 \\
yt&=M(x,z)(v^2+1)(1-xv-zv+xz)\end{array}\right.$$
Notons que par construction les classes de diviseurs 
$(\sqrt{a_i},0)-(-\sqrt{a_i},0)$
sont dans le noyau de l'endomorphisme de $J(C)$ dans $J(C')$ 
associ\'e \`a $\Gamma$, qui
contient donc le sous-groupe d'ordre $2^g$ de $J(C)[2]$ engendr\'e
par ces \'el\'ements. 

\medskip
Le th\'eor\`eme pour $g$ pair s'ensuit alors de la proposition suivante:

\medskip

\ni
{\sc Proposition.-} {\it Soit $\Gamma'\subset C'\times C$ la correspondance sym\'etrique de
$\Gamma$; $\Gamma'\circ \Gamma$ agit sur  $\Pic^0(C)$   par
$D\ms 2D$.}

\medskip

On montre sans difficult\'e, \`a partir des formules d\'efinissant $\Gamma$,   
que l'image d'un point $P=(X,Y)$ de $C$ par $\Gamma'\circ 
\Gamma$ est le diviseur $2P+P_1+w(P_1)$, o\`u $P_1$ est un point de $C$ d'abscisse
$-X$ et o\`u $w$ est l'involution hyperelliptique de $C$; l'action sur les classes de
diviseurs de degr\'e $0$ est donc la multiplication par $2$.

\subsection{Le cas $g$ impair}
Pour obtenir le th\'eor\`eme pour $g$ impair, il suffit, dans la construction
pr\'ec\'edente, de sp\'ecialiser $a_g$ en $0$; les courbes
$C$ et $C'$ sont alors de genre $g-1$, un calcul imm\'ediat donnant comme
\'equation de $C'$ celle indiqu\'ee dans le th\'eor\`eme.

\subsection{Dimension dans l'espace des modules}

1) Le cas $g=2$.

La courbe hyperelliptique g\'en\'erique de genre $2$ est du type $C$ ci-dessus; en effet, si $P_1,\ldots,P_6$
sont $6$ points g\'en\'eriques de la droite projective, il existe une unique involution $u$ telle que
$u(P_1)=P_2$ et $u(P_3)=P_4$; il existe alors une unique involution $w$, commutant
\`a $u$, telle que $w(P_5)=P_6$; dans un rep\`ere de la droite o\`u $u$ est 
donn\'ee par $x\ms -x$, $w$ est de la forme $x\ms t/x$, que l'on ram\`ene
\`a $x\ms 1/x$ par une homoth\'etie.

\ni
2) Le cas $g\geq 3$.

Deux courbes hyperelliptiques sont isomorphes si et seulement s'il existe
une homographie envoyant les points de Weierstrass de l'une sur ceux de l'autre.

Il suffit donc de prouver que, si $v,x_1,\ldots,x_g$  sont des points
g\'en\'eriques de $\P^1$, et si $h:\;x\ms (ax+b)/(cx+d)$ est une homographie 
telle que l'ensemble 
$A=\{h(v),h(1/v),h(x_1),\ldots,h(x_g),h(-x_1),\ldots, h(-x_g)\}$ 
est de la
forme $$\{w,1/w,y_1,\ldots,y_g,-y_1,\ldots,-y_g\},$$ $h$ 
est de la forme $x\ms \pm x$ ou $x\ms \pm 1/x$.

Soit 
$B=h^{-1}(\{y_1,-y_1,y_2,-y_2,u_3,-y_3\})$; $B$ est globalement
invariante par l'involution $h^{-1}uh$, o\`u $u$ est l'involution
$x\ms -x$.  

Or,  si $a_1,\ldots,a_6$ sont six \'el\'ements distincts d'un corps, il   
existe une involution permutant $a_{2i-1}$ et $a_{2i}$, $1\leq i\leq 3$
si et seulement si on a 

$$a_{{6}}a_{{5}}a_{{3}}+a_{{6}}a_{{5}}a_{{4}}+a_{{6}}a_{{2}}a_{{1}}+a_{{
5}}a_{{2}}a_{{1}}+a_{{1}}a_{{4}}a_{{3}}+a_{{2}}a_{{4}}a_{{3}}
$$ 
$$=a_{{6}}a_{{5}}a_{{1}}+a_{{6}}a_{{5}}a_{{2}}+a_{{6}}a_{{4}}a_{{3}}+a_{{
5}}a_{{4}}a_{{3}}+a_{{2}}a_{{1}}a_{{3}}+a_{{2}}a_{{1}}a_{{4}}
.$$

Par suite, tout \'el\'ement de $B$ est alg\'ebriquement d\'ependant des
autres; donc, si $b\in B$ est de la forme $\pm x_i$, on a $\{-x_i,x_i\}\subset 
B$, et si $b$ est \'egal \`a $v$ ou $1/v$, on a $\{v,1/v\}\subset B$. 
\`A une permutation de $\{1,\ldots,g\}$ pr\`es, $B$ est
donc de la forme 
$B_1=\{x_1,-x_1,x_2,-x_2,v,1/v\}$ ou
$B_2=\{x_1,-x_1,x_2,-x_2,x_3,-x_3\}.$ 

Comme prouv\'e plus haut, six points g\'en\'eriques de la droite
projective  s'\'ecrivent dans un rep\`ere convenable sous la forme
$\{x_1,-x_1,x_2,-x_2,v,1/v\}$, et donc, g\'en\'eriquement, il n'y a donc pas d'involution
conservant  $B_1$; donc $B$ est de la forme $B_2$ et $h(\{v,1/v\})=\{w,1/w\}$.
Or la courbe g\'en\'erique de genre
$2$ ayant comme groupe d'automorphismes $(\Z/2\Z)^2$ est justement
de la forme $y^2=(x^2-x_1^2)(x^2-x_2^2)(x^2-x_3^2)$, son groupe
d'automorphismes \'etant form\'e des quatre \'el\'ements
$(x,y)\ms (\pm x,\pm y)$. Il n'y a donc g\'en\'eriquement pas d'autre
involution que $x\ms -x$ conservant $B_2$; par suite, pour $1\leq i\leq 3$, 
$h^{-1}uh(x_i)=u(x_i)$; donc $h^{-1}uh=u$, et $h$ est une homographie
commutant \`a $u$, donc de la forme $x\ms ax$ ou $x\ms a/x$. Comme elle envoie
$\{v,1/v\}$ sur $\{w,1/w\}$, on a $a^2=1$, d'o\`u le r\'esultat. 

\subsection{Simplicit\'e de $J(C)$}
Pour $g=2$, la courbe $C$ \'etant la courbe g\'en\'erique de genre $2$, sa
jacobienne est absolument simple.

\ni
Pour $g=3$, on sp\'ecialise les ind\'etermin\'ees en prenant 
par exemple $v=2,a_1=1,a_2=3,a_3=4$; 
le polyn\^ome caract\'eristique du Frobenius en $13$ est
$${y}^{6}+2\,{y}^{5}+3\,{y}^{4}+44\,{y}^{3}+39\,{y}^{2}+338\,y+2197,$$
dont les racines sont 
$-(1+2i\cos\frac{5\pi}{7})(1+2i\cos \frac{3\pi}{7})(1+2i\cos\frac{\pi}{7})$
et ses conjugu\'ees, 
avec $i=\sqrt{-1}$;
le corps   qu'elles engendrent est le corps
$\Q(i,2\cos \frac{2\pi }{7})$, dont les racines de l'unit\'e
sont celles du corps
$L=\Q(i)$.   

Si la jacobienne n'\'etait pas absolument simple, 
il existerait un entier
$n$ tel que $y^n$ appartienne \`a $L$; $y^n$ serait \'egal, \`a une racine
de l'unit\'e pr\`es, \`a $(3\pm 2i)^n$; donc, \`a une racine de l'unit\'e
pr\`es, $y$ serait \'egal \`a un \'el\'ement de $L$, et $y$ serait dans $L$.

\medskip

\ni
Pour $g\geq 4$,  on raisonne par r\'ecurrence sur $g$: lorsqu'on sp\'ecialise $x_g$ en $0$, on trouve la courbe
$C$ de genre $g-1$ associ\'ee \`a $v,x_1,\ldots,x_{g-1}$; si $J(C)$ n'est pas
simple, elle est donc isog\`ene \`a $D\times E$, o\`u $D$ est  
absolument simple
de dimension $g-1$.

\medskip

Si l'on sp\'ecialise $v$ en $\sqrt{-1}$, la courbe $C$ 
admet comme automorphisme $(x,y)\ms (-x,y)$, et est rev\^etement de degr\'e $2$
des deux courbes d'\'equation
$y^2=(x+1)(x-x_1^2)\ldots (x-x_g^2)$ 
et 
$y^2=x(x+1)((x-x_1^2)\ldots (x-x_g^2)$), de genre $g/2$ si $g$ pair et
de genre $(g-1)/2$ et $(g+1)/2$ sinon; $J(C)$ est donc isog\`ene au
produit de leurs jacobiennes, qui sont g\'en\'eriquement absolument simples.
Ceci contredit le fait que $J(C)$ soit isog\`ene \`a $D\times E$;
par suite, d\`es que $g\geq 4$, $J(C)$ est absolument simple.

\bigskip

{\sc R\'ef\'erences}

[1] J.-B. Bost and J.-F. Mestre. {\it Moyenne arithmético-géometrique et 
périodes de courbes de genre 1 et 2}. Gaz. Math. Soc. France $38$, 1988, 36-64. 

\medskip

[2] F. Richelot. {\it Essai sur une m\'ethode g\'en\'erale pour
d\'eterminer la valeur des int\'egrales ultra-elliptiques,
fond\'ee sur des transformations remarquables de ces transcendantes},
C.R. Acad. Sci. Paris 2, 1836, 622-627.

\medskip 

[3] F. Richelot . {\it De transformatione integralium Abelianorum 
primiordinis commentation}, J. Reine Angew. Math. $16$, $221-341$.

\end{document}